\documentclass{aip-cp}

\usepackage[numbers]{natbib}
\usepackage{rotating}
\usepackage{graphicx}

\usepackage{amsmath}

\newcommand{\x}{\mathbf{x}}
\newcommand{\q}{\mathbf{q}}
\newcommand{\CFL}{\textnormal{CFL}}

\newcommand{\R}{\mathbb{R}}
\newcommand{\f}{\mathbf{f}} 

\newcommand{\Q}{\mathbf{Q}}

\newcommand{\F}{\mathbf{F}}
\newcommand{\B}{\mathbf{B}}

\newcommand{\N}{\mathbb{N}}
\newcommand{\n}{\mathbf{n}}

\renewcommand{\vec}[1]{\mathbf{#1}}

\newcommand{\diff}{{\mathrm{d}}}

\begin{document}

% Title portion
\title{Second Order Finite Volume Scheme for Shallow Water Equations on Manifolds}
%\title{The Title Goes Here with Each Initial Letter Capitalized}

\author[aff1]{Michele Giuliano Carlino}
\eaddress{michele-giuliano.carlino@inria.fr}
\author[aff1]{Elena Gaburro\corref{cor1}}
%\eaddress{anotherauthor@thisaddress.yyy}

\affil[aff1]{Inria, Univ. Bordeaux, Bordeaux INP, UMR 5251, 200 Avenue de la Vielle Tour, 33405 Talence cedex, France.}
%\affil[aff2]{Additional affiliations should be indicated by superscript numbers 2, 3, etc. as shown above.}
%\affil[aff3]{You would list an author's second affiliation here.}
\corresp[cor1]{Corresponding author: elena.gaburro@inria.fr}

\maketitle

\begin{abstract}
In this work we propose a second-order accurate scheme for shallow water equations in general covariant coordinates over manifolds. In particular, the covariant parametrization in general covariant coordinates is induced by the metric tensor associated to the manifold. The model is then re-written in a hyperbolic form with a tuple of conserved variables composed both of the evolving physical quantities and the metric coefficients. This formulation allows the numerical scheme to i) automatically compute the curvature of the manifold as long as the physical variables are evolved and ii) numerically study complex physical domains over simple computational domains.
\end{abstract}

% Head 1
\section{INTRODUCTION}
In this work, we develop a second-order accurate scheme \cite{gaburro2017direct, bergmann2022second} for the Shallow Water (SW) system in general covariant coordinates over a manifold. SW equations find application in geophysics for studying and modelling river and ocean flows, meteorology, avalanches and tsunamis. For all these applications, the scale dimension is different. In particular, the definition of the gravity field changes case by case. Once an equipotential surface is detected and assumed to be the manifold, a relative metric tensor is deduced. As a consequence, the differential operator of the system is affected by the curvature of the manifold through the metric. We propose a mathematical model where the metric is given as an initial datum. Then, the hyperbolic re-formulation of the original SW system automatically computes the evolution in time of the related physical variables as well as the geometric characteristics associated to the manifold. According to the shallow water hypothesis, the movement of fluid particles along the direction of gravity is neglected. In general, the gravitational force is considered oriented vertically, also for large distances, namely when the gravity field is not uniformly parallel to a straight line. In \cite{bachini2020geometrically}, the authors numerically proved that this approximation strongly affects the final solution, which results to be far from the one obtained with the physical gravity field, even for almost constant bathymetries. Thus, this work addresses a formulation of the problem allowing to easily account for a physical solution of the phenomenon. By following a previous work \cite{gaburro2021well} for general relativity and as largely studied in \cite{carlino2022well}, the terms related to the spatial derivatives in general covariant coordinates are not considered as source terms of the system but properly defined in the nonconservative components of the arising hyperbolic re-formulation. Thus, the state variable of the problem is composed of physical quantities (i.e. the fluid depth and velocity) and metric quantities. This formulation has a double advantage: it allows to evolve the system with an automatic computation of metric-related variables and makes the system frame-free with respect to the general covariant coordinates used for describing the manifold. Consequently, the numerical system is solved over a computational (reference) domain simpler w.r.t. the physical domain. The present work can also be seen as an extension to general covariant coordinates of the work \cite{castro2017well} valid for spherical coordinates.

This contribution is organized as follows. First, the free-frame re-formulation of SW system is introduced, followed by the minimal description of the second order Finite Volume (FV) adopted strategy. A numerical validation is showed before the conclusions. Along the whole contribution, the considered spatial dimension is 2. Moreover, the Einstein notation for the sum on repeated indexes is used.

\section{FRAME-FREE FORMULATION OF SHALLOW WATER SYSTEM ON MANIFOLD} \label{sec.formulation}
The system of equations to be solved is written in a hyperbolic form of conservation law as
\begin{equation} \label{eq.hyp}
	\partial_t \Q + \nabla \cdot \F(\Q) + \B(\Q) \cdot \nabla \Q = \mathbf{0}, \quad \x = (x^1, x^2) \in \Omega, \quad t \in [0, T].
\end{equation}
In the above expression, $\x$ is the vector space in the computational domain $\Omega \subset \R^2$ and $t$ is time in interval $[0, T]$. The conserved variables are collected in the state $\Q$ belonging to the space of admissible states $\Omega_\Q \in \R^m$. Terms $\F(\Q) = (\f_1(\Q), \f_2(\Q))$ and $\B(\Q) = (\B_1(\Q), \B_2(\Q))$ are the nonlinear flux and the nonconservative components, respectively, with $\f_i \in \R^m$ and $\B_i \in \R^{m \times m}$ ($i = 1,2$). System (\ref{eq.hyp}) is hyperbolic if, for any direction $\n = (n^1, n^2)$, matrix $\mathbf{A} = (\partial \F / \partial \Q + \B ) \cdot \n$ has only real eigenvalues whose corresponding eigenvectors are linearly independent. 

The original SW system on an arbitrary manifold $\mathcal{M}$ in general covariant coordinates reads \cite{baldauf2020discontinuous, carlino2022well}
\begin{equation} \label{eq.baldauf}
	\begin{aligned}
		\frac{\partial h}{\partial t} + \frac{\partial m^j}{\partial x^j} + \Gamma^j_{jk} m^k &= 0,\\
		\frac{\partial m^i}{\partial t} + \frac{\partial T^{ij}}{\partial x^j} + \Gamma^i_{jk} T^{kj} + \Gamma^j_{jk} T^{ik} &= S^i; \quad i = 1,2.
	\end{aligned}
\end{equation}
The first and second lines in \eqref{eq.baldauf} are the balance for the mass and the momentum, respectively. They describe the evolution of the fluid depth $h$ and the mass flux $m^i = h u^i$, where $u^i$ is the $i$-th component of the fluid velocity. In this model, $h$ is measured as the length of the gravitational line from the free surface of the fluid to the bottom. Accordingly, the direction of the velocity components $(u^1, u^2)$ are given by the frame $(x^1, x^2)$. In the momentum equation the flux is led by the stress tensor $T^{ij} = m^i m^j/h + 0.5 g h^2 \gamma^{ij}$, where $g$ is the gravitational constant and $\gamma^{ij}$ is the metric tensor in contravariant representation. Moreover, the source term is defined as $S^i = -g h \gamma^{ij} \partial_{\partial x^j} b$, with $b$ the bathymetry. In the source term, any force due to Coriolis' effect, bottom-friction phenomena or viscosity is not considered. Due to the curvature of the manifold $\mathcal{M}$, the differential operator of system \eqref{eq.baldauf} is also defined by the Christoffel symbols $\Gamma^i_{jk} = 0.5 \gamma^{im} ( \partial_{\partial x^j} \gamma_{km} + \partial_{\partial x^k} \gamma_{jm} - \partial_{\partial x^m} \gamma_{jk} )$, with $\gamma_{ij}$ the covariant representation of the metric tensor.

By following \cite{carlino2022well}, after exploiting the contraction properties of the Christoffel symbols, properly applying the chain rule on physical flux terms $\partial_{x^j} m^j$ and $\partial_{x^j} T^{ij}$ and lowering the indexes of the contravariant components of the metric (due to the relation $[\gamma_{ij}]^{-1} = [\gamma^{ij}]$), system \eqref{eq.baldauf} writes in form \eqref{eq.hyp} with the state $\Q = [ h, m^1, m^2, b, \gamma_{11}, \gamma_{12}, \gamma_{22} ]^T$; the flux components $\f_1 = [m^1, (m^1)^2/h, m^1 m^2/h, \mathbf{0}_4^T]^T$ and $\f_2 = [m^2, m^1 m^2/h, (m^2)^2/h, \mathbf{0}_4^T]^T$; and the nonconservative matrices
\begin{equation} \label{eq.NCP}
	\begin{aligned}
		\B_1 &= \frac{1}{\gamma} \left[ \begin{array}{c}
			\begin{matrix}
				0 & 0 & 0 & 0 & \frac{1}{2} m^1 \gamma_{22} & -m^1 \gamma_{12} & \frac{1}{2} m^1 \gamma_{11} \\[3pt]
				g h \gamma_{22} & 0 & 0 & g h \gamma_{22} & \frac{(m^1)^2}{h} \gamma_{22} & -2 \frac{(m^1)^2}{h} \gamma_{12} & -\frac{1}{2 h} [2 m^1 m^2 \gamma_{12} + \alpha_{11}] \\[3pt] 
				- g h \gamma_{12} & 0 & 0 & - g h \gamma_{12} & \frac{m^1}{2 h} (-m^1 \gamma_{12} + m^2 \gamma_{22}) & \frac{m^1}{h} (m^1 \gamma_{11} - m^2 \gamma_{12}) & \frac{m^2}{2 h} (3 m^1 \gamma_{11} + m^2 \gamma_{12})
			\end{matrix} \\[3pt] \hline
			\mathbf{O}_{4 \times 7}
		\end{array} \right], \\
		\B_2 &= \frac{1}{\gamma} \left[ \begin{array}{c}
			\begin{matrix}
				0 & 0 & 0 & 0 & \frac{1}{2} m^2 \gamma_{22} & -m^2 \gamma_{12} & \frac{1}{2} m^2 \gamma_{11} \\[3pt]
				- g h \gamma_{12} & 0 & 0 & - g h \gamma_{12} & \frac{m^1}{2 h} (3 m^2 \gamma_{11} + m^1 \gamma_{12}) & \frac{m^2}{h} (m^2 \gamma_{11} - m^1 \gamma_{12}) & \frac{m^2}{2 h} (-m^2 \gamma_{12} + m^1 \gamma_{22}) \\[3pt] 
				g h \gamma_{11} & 0 & 0 & g h \gamma_{11} & -\frac{1}{2 h} [2 m^1 m^2 \gamma_{12} - \alpha_{22}] & -2 \frac{(m^2)^2}{h} \gamma_{12} & \frac{(m^2)^2}{h} \gamma_{22}
			\end{matrix} \\[3pt] \hline
			\mathbf{O}_{4 \times 7}
		\end{array} \right],
	\end{aligned}
\end{equation} 
with $\gamma = \det[\gamma_{ij}] = \gamma_{11} \gamma_{22} - \gamma_{12}^2$ and $\alpha_{ij} = -(m^1)^2 \gamma_{ii} + (m^2)^2 \gamma_{jj}$. In the above expressions, $\mathbf{0}_4$ and $\mathbf{O}_{4 \times 7}$ are the null vector and the null matrix in $\R^4$ and $\R^{4 \times 7}$, respectively. The state vector is composed of four physical variables $(h, m^1, m^2, b)$, i.e. the fluid depth $h$ and the mass fluxes $m^i$, and three metric components in covariant representation $(\gamma_{11}, \gamma_{12}, \gamma_{22})$. In particular, the metric components are three because of the symmetry of the metric tensor.

\section{SECOND ORDER FINITE VOLUME MUSCL-HANCOCK SCHEME} \label{sec.scheme}
A second order FV MUSCL-Hancock-type approach \cite{van1979towards, toro2013riemann} is employed. For more details, see \cite{carlino2022well} and its bibliography. Let $\mathcal{T}_N$, with $N \in \N$, be a Voronoi-type tessellation \cite{boscheri2022continuous, gaburro2020high} for the computational domain $\Omega$. In particular, the spatial discretization defines a set of mutually disjoint cells $\Omega_k \in \mathcal{T}_N$, with center of mass $\x_k$ for $k = 1, \dots, N$, whose union fully covers the domain. To any cell $\Omega_k$, a stencil $\mathcal{S}_k$ of cells sharing at least one vertex is defined. Moreover, the edge of two contiguous cells $\Omega_k$ and $\Omega_l$ is denoted by $e_{kl}$. The set of edges of a cell $\Omega_k$ is denoted by $\mathcal{E}_k$. Time set is split into segments $(t^n, t^{n+1})$ such that $t^{n+1} = t^n + \Delta t$, with $\Delta t$ the time step. Finally, the cell average of the state at time $t^n$ is defined as $\Q_k^n = |\Omega_k|^{-1} \int_{\Omega_k} \Q(\x, t^n) \, \diff \x$.

The first step of the scheme consists in a local space-time reconstruction of the numerical solution, known at time $t^n$ in the form of cell averages. The reconstruction reads
\begin{equation} \label{eq.reconstruction}
	\q_k(\x, t) = \mathbf{w}_k(\x) + \partial_t \Q_k \, (t - t^n) = \Q_k^n + \nabla \Q_k \, (\x - \x_k) + \partial_t \Q_k \, (t - t^n); \quad \x \in \Omega_k, \quad t \in (t^n, t^{n+1}), 
\end{equation}
where spatial and time polynomial coefficients $\nabla \Q_k$ and $\partial_t \Q_k$ in $\R^m$ have to be found. For the spatial coefficient $\nabla \Q_k$, the integral conservation of $\mathbf{w}$ over the stencil $\mathcal{S}_k$ is imposed, by applying a slope limiter \cite{toro2013riemann, van1979towards}. Concerning the time coefficient $\partial_t \Q_k$, it is found by integrating in discrete sense system \eqref{eq.hyp} along the boundary $\partial \Omega_k = \bigcup_{e_{kl} \in \mathcal{E}_k} e_{kl}$ of cell $\Omega_k$. Once reconstruction \eqref{eq.reconstruction} is computed, the cell average evolution from time $t^n$ to $t^{n+1}$ is found via the FV scheme
\begin{equation} \label{eq.muscl-hancock}
	\Q_k^{n+1} = \Q_k^n - \frac{\Delta t}{|\Omega_k|} \sum_{e_{kl} \in \mathcal{E}_k} |e_{kl}| \boldsymbol{\mathcal{F}}_{e_{kl}}(\q_{e_{kl}}^-, \q_{e_{kl}}^+) \cdot \n_{e_{kl}} - \frac{\Delta t}{| \Omega_k |} \sum_{e_{kl} \in \mathcal{E}_k} \mathcal{D}_{e_{kl}}(\q_{e_{kl}}^-, \q_{e_{kl}}^+) - \Delta t \B(\q_k^{n + \frac{1}{2}}) \cdot \nabla \Q_k,
\end{equation}
where $\boldsymbol{\mathcal{F}}_{e_{kl}}$ and $\mathcal{D}_{e_{kl}}$ are the numerical flux approximation and nonconservative jump at the interface, respectively. Their approximation exploits the local reconstructions \eqref{eq.reconstruction} from the right ($\q^+$) and from the left ($\q^-$) with respect to any edge $e_{kl}$.  The flux is approximated via a classical Rusanov approach. Concerning the nonconservative jumps at the interface, they are treated by integrating along the Lipshitz-continuous segment path  $\mathbf{\Psi}(\q^-, \q^+; \tau) = \q^- + \tau (\q^+ - \q^-)$ from state $\q^-$ to $\q^+$, with $\tau \in [0,1]$, following \cite{dal1995definition,pares2006numerical}. Finally, time step $\Delta t$ follows a standard $\CFL$ condition.

\section{NUMERICAL VALIDATION} \label{sec.validation}

\begin{figure}[!bp]
	\includegraphics[width=0.4\textwidth]{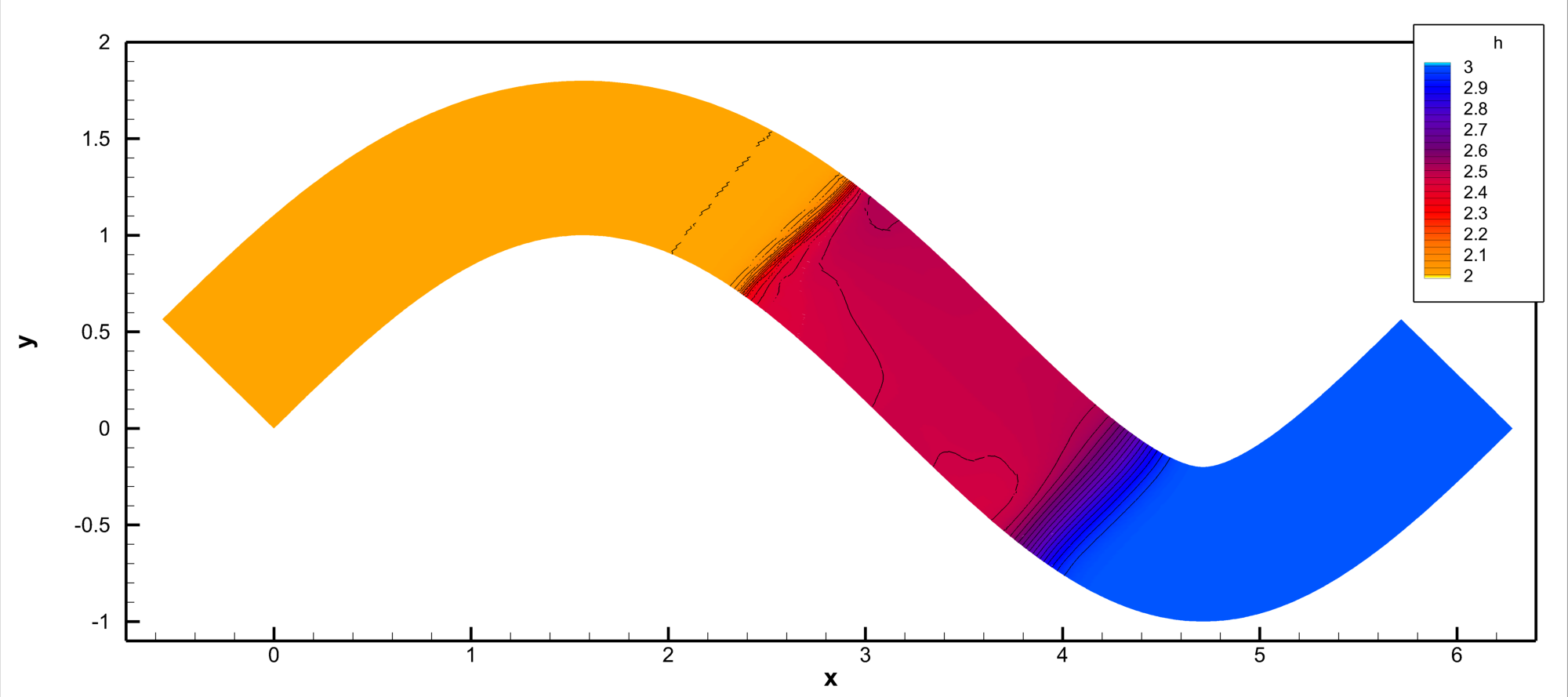}% 
	\includegraphics[width=0.4\textwidth]{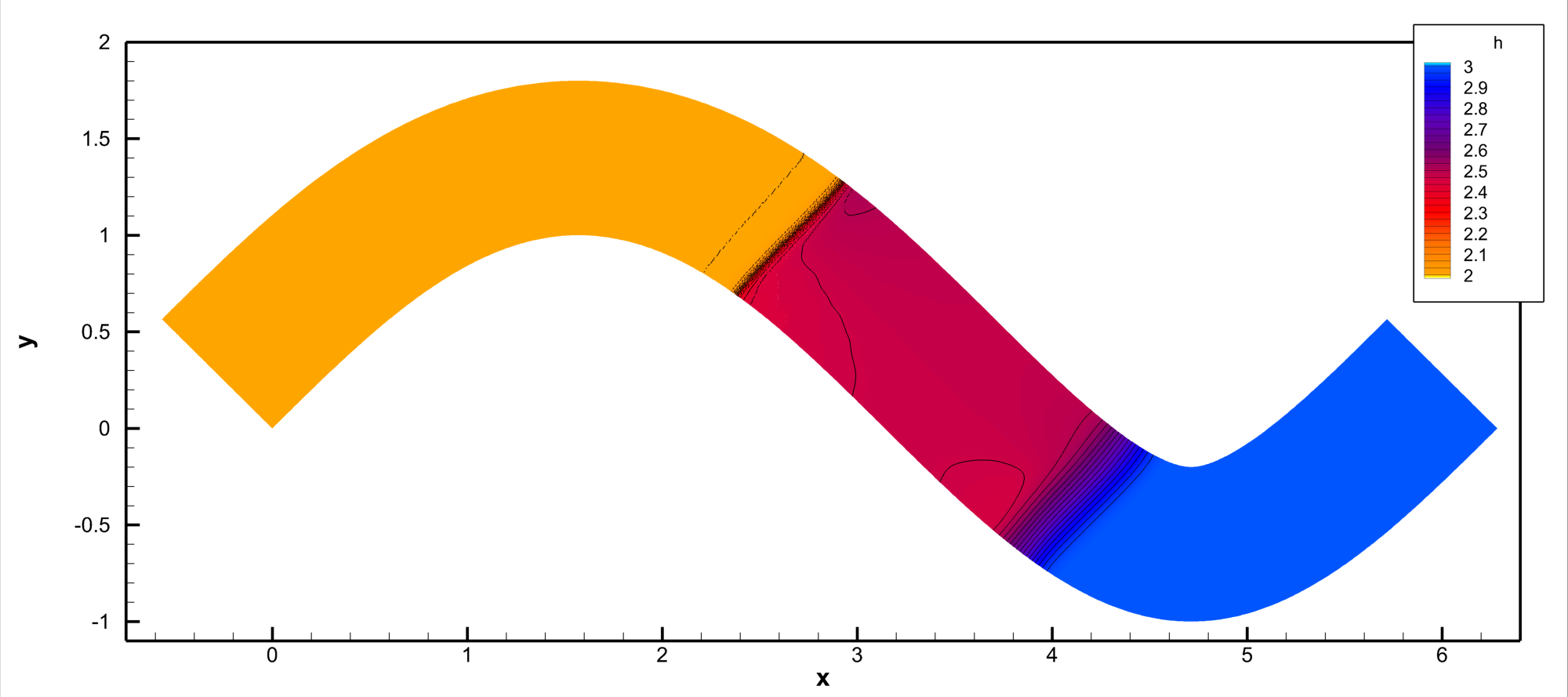}
	\label{fig.coarse}
	\caption{Solutions at final time $T = 0.2$ over a mesh of size 1.61E-2 from classical SW system in Cartesian coordinates (left) and from our novel approach in general covariant coordinates (right). }
\end{figure}
\begin{figure}[!bp]
	\includegraphics[width=0.4\textwidth]{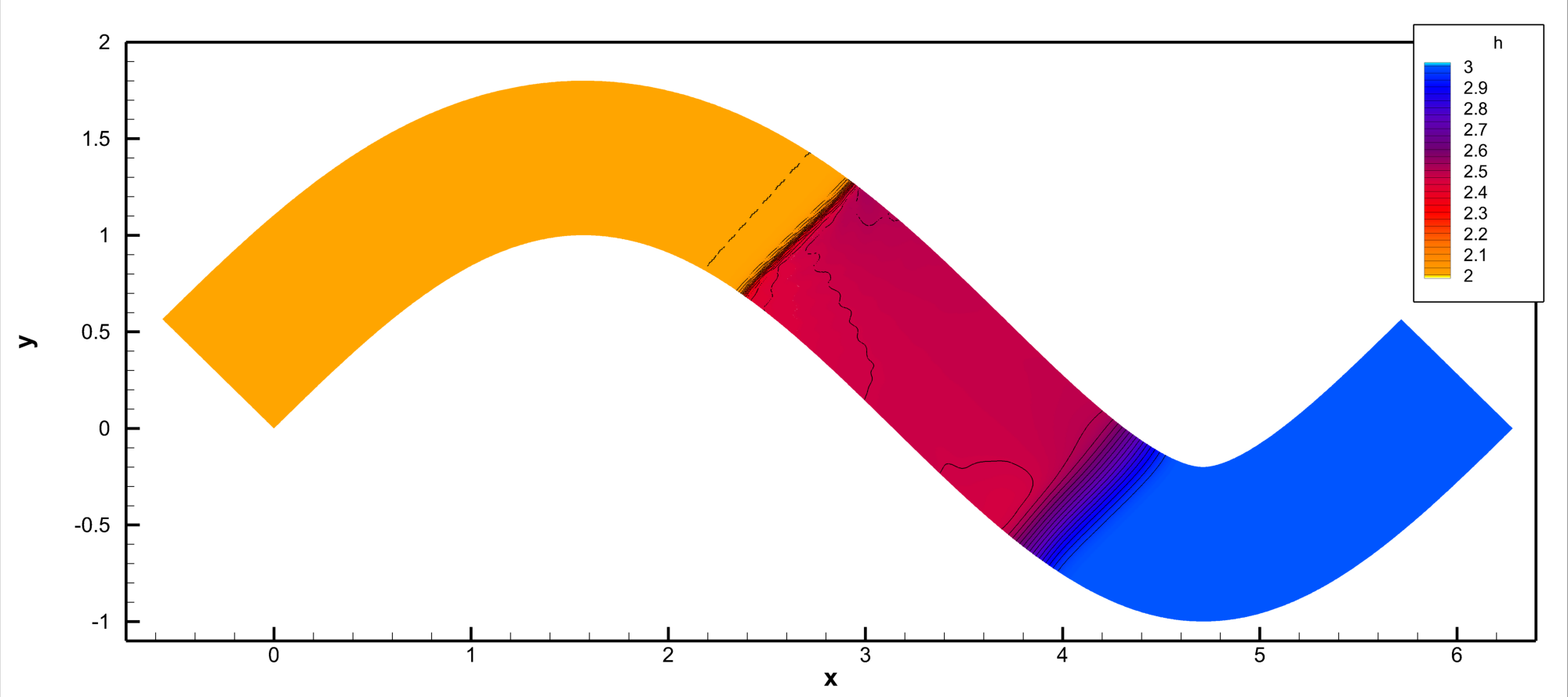}%
	\includegraphics[width=0.4\textwidth]{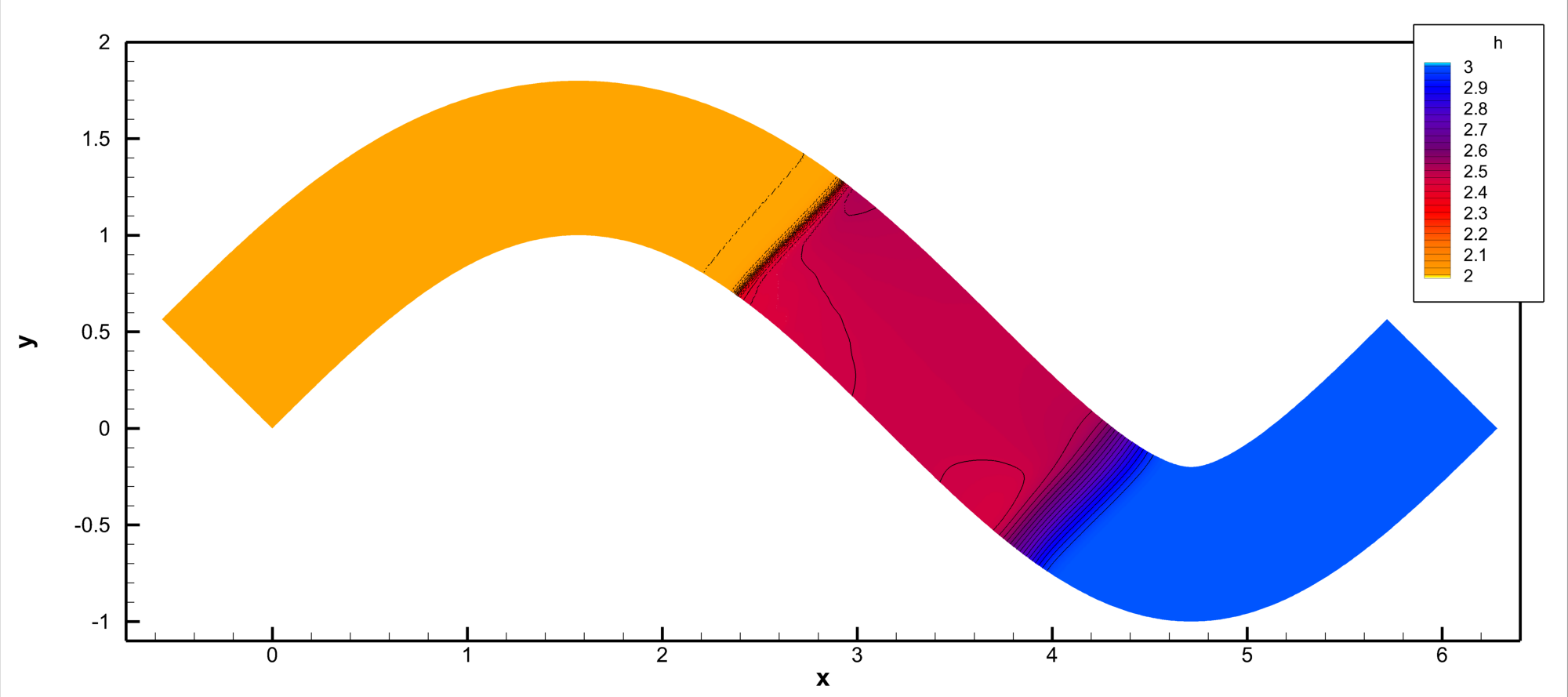}
	\label{fig.fine}
	\caption{Solutions at final time $T = 0.2$ over a mesh of size 4.16E-3 from classical SW system in Cartesian coordinates (left) and from our novel approach in general covariant coordinates  (right). }
\end{figure}

Let the computational domain $\Omega = [0, 2\pi] \times [0, 0.8]$ of general covariant coordinates $(x^1, x^2)$ be transformed in the S-shaped manifold $\mathcal{M}$ of Cartesian coordinates $(x,y)$ through the map
\begin{equation} \label{eq.map}
	x = x^1 - \frac{x^2 \cos x^1}{\sqrt{1 + \cos^2 x^1}}; \quad y = \sin x^1 + \frac{x^2}{\sqrt{1 + \cos^2 x^1}}.
\end{equation}  
Map \eqref{eq.map} defines a metric tensor whose entries are $(\gamma_{11}, \gamma_{12}, \gamma_{22}) = ([ (1 + \cos^2 x^1)^{\frac{3}{2}} + x^2 \sin x^1 ]^2/(1 + \cos^2 x^1)^2, 0, 1)$. The bathymetry $b \equiv 0$ is flat and the gravity is normally oriented w.r.t. the manifold $\mathcal{M}$ with modulus $g = 9.81$. We study the evolution of a Riemann problem with initial state defined by zero velocity and fluid depth $h_0(x,y) = 2\chi_{\{ y \geq x - \pi \}}(x,y) + 3\chi_{\{ y < x - \pi \}}(x,y)$, with $\chi$ the indicator function. Wall boundary conditions are imposed.
We compare the solutions at time $T = 0.2$ obtained with the classical SW system in the S-shaped physical domain, i.e. the manifold, and with our novel approach on the simple computational domain $\Omega$ by using the described metric quantities. Figures~\ref{fig.coarse} and~\ref{fig.fine} depict the comparison on meshes of size 1.61E-2 and 4.16E-3, respectively. The solutions obtained with the two approaches are in perfect agreement and two mesh sizes are reported to show the grid convergence.

\section{CONCLUSIONS} \label{sec.conclusions}
In this contribution we present a second-order FV scheme for SW equations over a manifold in general covariant coordinates. In particular, a frame-free re-formulation of the original system is proposed. The arising hyperbolic system allows to evolve the involved physical quantities as well as to easily compute the geometrical variables associated to the manifold. % (e.g. the Christoffel symbols). 
The state variable of the new problem collects the physical quantities and the metric tensor coefficients, thus the nonconservative associated components allow to take into account the curvature of the manifold. As a direct consequence, the solution can be computed on a reference computational domain simpler than the physical domain. In particular, the showed numerical test case accounts for a S-shape physical domain, but the numerical solution is found over a simple rectangle.
Further validation of the proposed model is presented in \cite{carlino2022well} where the results obtained with our formulation have been i) compared with SW written in Cartesian, polar and elliptical coordinates, and ii) checked for convergence and iii) for shock capturing capabilities.
Future works will exploit the developed technology in terms of treating differential geometry terms to improve the numerical simulation of systems involving general relativity \cite{gaburro2021well} for which the metric of the space--time is naturally playing a role and is also evolving in time.% 

% Acknowledgement
\section{ACKNOWLEDGMENTS}

E.~Gaburro is member of the CARDAMOM team at the Inria center of the University of Bordeaux and 
gratefully acknowledges the support received from the European Union’s Horizon 2020 Research and Innovation Programme under the Marie Sk\l{}odowska-Curie Individual Fellowship \textit{SuPerMan}, grant agreement No. 101025563. 

% References

%\nocite{*}
\bibliographystyle{plain}%
\bibliography{references}%

\end{document}